# Research on Real-time Operational Status Evaluation Technology of Automobile Based on Information Data Fusion Algorithm

*Zuo Yanhong, Xia Shilong, Zhou Chao, Yang Kun*

(School of Mechanical and Electrical Engineering, Anhui Jianzhu University, Hefei 230601,China)

**Abstract**：Automobiles have become the main means of transportation for human beings, and their failures in the process of operation are directly related to the life and property safety of drivers. Therefore, real-time operational status evaluation technology have become urgent problems in the current academic world. The premise of real-time operational status evaluation technology of automobiles is to obtain high-quality information data in real time, but the automobile operating environment is complex and changeable, resulting in the measured information data under the influence of multiple factors, such as equipment performance and signal interference. There is an unpredictable measurement error, which greatly affects the reliability of operational status evaluation systems. In this paper, on the basis of studying the structure and operation characteristics of automobiles, we design a method that can be used for real-time operational status evaluation of automobiles; through the study of fractional-order calculus theory, we establish a mathematical model of information data fusion based on fractional-order differential operators. By providing high-quality information data to automobile real-time operational status evaluation systems, real-time operational status evaluation technology for automobile faults can be realized. The feasibility and effectiveness of the method were verified through an experiment applying the technology in automobile real-time operational status evaluation systems. The research results are of great significance for promoting the development of the automobile industry and ensuring the safety of drivers' lives and property.

## 0 Introduction

With the development of society, the car has become the main means of transportation for human beings, from the moment the world's first car was born, according to incomplete statistics, the global automobile production has been more than 80 million, it can be seen that the car has become an indispensable part of people's lives. The development process of the automobile industry is actually the development process of science and technology, and with the development of science and technology along the trajectory of manual operation, mechanization, automation, and intelligence, the automobile has also developed from the steam era, internal combustion engine era, electronics and intelligence, and new energy power era. After comparing the performance and technical content of locomotives at different times, we can easily find that: with the advance of the times, the structure, function and performance of automobiles are undergoing qualitative transformation. Currently, all systems and components of automobiles are beginning to adopt electronic control and information processing. Automobiles have become a product of the multidisciplinary integration of engineering materials, mechanical processing, Internet of Things transmission, computer control, and intelligent manufacturing, which also creates difficulties for online inspection and fault diagnosis of automobiles.

Automobile real-time operational status evaluation technology play a vital role in the process of automobile use, and can provide technical support for automobile maintenance and repair by accurately evaluating the running state of automobiles. Currently, the assessment of the operational status of the vehicle is limited to the fault diagnosis level, where the comparison of the measured values of the operating parameters of the source of the fault with the standard values determines whether the source of the fault is in a faulty state or not. Methods of fault diagnosis can be summarized as follows [1,2].

(1) Manual empirical methods

Manual experience diagnostic method maintenance personnel can see, listen, smell, and other methods to determine where the vehicle has a problem, with the help of their own maintenance experience, as well as some simple tools for analysis. This method has been very practical until now, with low cost, simplicity, and convenience. At the same time, there are also disadvantages such as low efficiency and reliability, which are too dependent on the technical level of the technicians. 。

(2) Instrumental diagnostic methods

The instrumental diagnostic method refers to the maintenance personnel using testing instruments to determine the relevant parameter values, through the comparison of the measured value with the standard value, online diagnosis of automobile operation status. With improvements in the performance of equipment and instruments, the instrument method has experienced the development of a traditional simple instrument diagnosis method for precision instrument diagnosis. Compared with the manual empirical method, this method has the advantages of high working efficiency and accurate information data, but with low cost effectiveness, the reliability of fault diagnosis still relies on the technician's insufficiency.

(3) On-board system diagnostic methods

The on-board diagnostic system analyzes the data of automobile parts using the vehicle controller, and then displays it on the vehicle screen, or shows the parts and time of fault occurrence on the handheld diagnostic instrument. The on-board diagnostic system can carry out real-time monitoring of the vehicle, issue fault warning information, has a strong fault-tolerant function, and reduces the burden on maintenance personnel. However, the on-board diagnostic system can only diagnose some existing specific faults, the scope is limited, it cannot meet the actual demand, and often results in misjudgment. When a control unit fails, it cannot diagnose itself and needs to be assisted by other diagnostic techniques.

(4) Artificial Intelligence Diagnostics

Owing to the many shortcomings of fault diagnosis, with the development of science and technology, artificial intelligence has become a hot research topic. Automotive artificial intelligence diagnostic technology is an emerging comprehensive technology that integrates the theories and techniques of various disciplines such as computer science, artificial intelligence, pattern recognition, automatic control theory, information analysis and processing, and electronic technology [3]. Thus far, automotive intelligent fault diagnosis methods are mainly divided into two categories: one is based on mathematical model fault diagnosis methods, and the other is artificial intelligence diagnosis. The main methods include artificial neural networks [4], Bayesian networks [5], expert systems [6], and other diagnostic methods. The main research results are: Mao Ziheng applied an extreme sample entropy (ESE)-based application to the diagnosis of minor short-circuit faults in automotive lithium-ion batteries [7]. Wang Jianping diagnosed automotive permanent magnet synchronous motors (PMSMs) and permanent magnet synchronous motors (PMSMs) by combining the empirical mode decomposition-symmetric point pattern (EMD-SDP) image features and an improved DenseNet convolutional neural network. PMSM (permanent magnet synchronous motor) fault diagnosis function [8]. Li Xiwei proposed a spectral self-focusing fault diagnosis method for automotive transmissions under gear shifting conditions to realize the fault diagnosis function of automotive transmissions under gear shifting conditions [9]. Salah and Alia combined parameter estimation and machine-learning classification algorithms to provide a complete model mapping of the health and fault states of a machine to implement a fault diagnosis technique for the main motors of a vehicle [10]. Jawalkar utilized machine learning algorithms to implement fault detection and diagnostic functions in automotive control systems [11]. Melis and Tommaso utilized silicon Melis, and Tommaso utilized the basic principle of photoemission of electronic components to implement the analog circuit fault diagnosis function in a team car by inserting basic devices that can be used to check the state of the circuit for each basic signal photoemissioning return circuit [12]. Nair, S introduced a pioneering condition monitoring methodology rooted in vibration analysis, in which a set of auto-regressive moving averages (ARMA), histograms, and statistical features are carefully extracted from vibration signals acquired during operation, to develop a technique that can be used specifically for detecting faults in clutch systems [13]. The Supervised Feature Extraction Method utilizes the global structure of the data extracted from a large number of unlabeled samples to adjust Fisher's criterion to achieve the function of detect the faults generated in the automobile assembly process [14]. Liang, C Generalized composite multiscale diversity entropy and its application for fault diagnosis of rolling bearings in automotive production lines [15].

According to the analysis of the existing automotive online detection and fault diagnosis methods, it can be seen that they have the following disadvantages in engineering applications:

① Each method is only adapted to a certain type of automobile fault diagnosis field and cannot be applied to other types of fault diagnosis fields; therefore, there are limitations in the field of application.

② Failure to take into account the complexity of the working environment in the process of automobile operation; ignoring the working environment affects the performance of the equipment as an important factor, so there is a lack of reliability drawbacks.

③ Their research content is biased in the field of monitoring and fault diagnosis of the information during the operation of the car, and they do not have fault prediction engineering, which leads to the lagging maintenance of the car.

More importantly: the above methods are only limited to the application level of automobile fault diagnosis, and have not yet realized the early warning function of automobile faults, therefore, they cannot avoid the occurrence of its faults through the method of real-time judgment of automobile's operation status, and through the method of fault early warning.

As a result, no ideal real-time operational status evaluation technology for automobiles has been found until now. In this paper, we intend to combine the research results in the literature [16] and attempt to combine computer technology and High-precision information data detection to realize accurate monitoring and fault diagnosis technology for all types of automobile accident information.

# 1. Real-time operational status evaluation technology for automobiles

## 1.1 Working Principle

1.1.1 Individual working state judgment and early warning

Assuming that there exists a test object whose measured value $J$ changes regularly according to the function $J(x)$, where $x$ is the influence factor of the measured value and $x \in [a,b]$, the interval of change in the measured value $F$ is $[J_a, J_b]$. When the measured value is $J_i$, the parameter value $x_i$ of its influence factor $x$ can be calculated according to formula $F(x)$, and according to the above functional relationship, the margin of the influence factor $x$ from the extreme value point $J_b$ can be derived. To ensure that detection object A is in the normal working state, set the safety coefficient $k$ and $k \in (0,1)$ of the measured value. When, its measured value $J = k*b$ will be in the warning state.

2.1.2 Judgment and warning of system working state

Suppose there exists a working system, which has n mutually independent detection objects. Each detection object measurement value has its own change rule $J(x)$ and $x \in [a,b]$, then each detection object value interval $J=[J(b)-J(a)]$. n detection objects are running normally at the same time, and at a certain point in time $k$, the measurement value of each measurement object $J_i$, where $i \in [1,n]$, $k \in [a_i,b_i]$, and combined with the equation $J(x)$, the interval $T_i$ of the influence factors of each detection object from its limit value can be obtained as

$$\begin{cases} T_1 = b_1 - k \\ T_2 = b_2 - k \\ T_3 = b_3 - k \\ \dots \\ T_n = b_n - k \end{cases} \quad (1)$$

Because the system consists of n detection objects working in parallel, as long as the measurement value of one of these objects is not within the normal range, the entire system is in an abnormal working state. Therefore, the time of failure of the system state is the minimum value $T_{min}$ in (1), and the corresponding detection object is the source of its failure. The best way to avoid the occurrence of faults is the fault warning, in the fault is about to occur when the system will issue a warning signal to remind the user to timely maintenance of equipment to avoid the occurrence of faults. The failure warning function can be realized by setting the detection object warning measurement value $J_y$ method and setting the safety coefficient $k$ of each measurement object measurement value, $k \in (0,1)$; then, the detection object Q measurement value $J_y = k * J_e$, where $J_e$ is the rated lifeand the entire system will be in a warning state.

## 1.2 Troubleshooting and Early Warning Methods

Automobiles are composed of multiple systems, such as power, transmission, and air conditioning systems, and each system is composed of multiple components. Each component has its own life cycle $T$, and the life cycles $T$ of different components are different. Therefore, in the normal operation of a car, each component works at the same time, and the different components are in different positions in their life cycles. Assuming that the car has n important parts, each part of the life cycle of $T$, the main parameters affecting its life cycle, and the time t of the functional expression for $f(t)$, and $0 < t \leqslant T$. If in the normal operation of the car at a certain point in time $t_i$ through the on-board information data detection system to measure the main parameters of the component value $f_i$, the function $f(x)$ can be used to derive the part located in the life cycle of the time point $t_i$, the part located in the life cycle of the time point $t_i$, and the car located in its life cycle of the time point $f(x)$. The life cycle of the time point $t_i$, then the time point of the failure of the component from the measurement of the time point after $T_i = T - t_i$.

The time of each important component from its life cycle is $T_1, T_2 \dots T_n$. For an entire vehicle, as long as one of the major components fails, the entire vehicle will be in a failure state. Therefore, the vehicle failure time points for $T_1, T_2\dots$ in the smallest value of $T_{min}$, the source of the failure is the corresponding part. To avoid the loss of life and property when a car fails, the automobile control system should have a failure warning function. Setting the failure warning time point reminds the car owner that the car is close to the failure state, so that he can maintain the car in time. By setting the warning time coefficient of each important part $k$, and $k \in (0,1)$, assuming that the life cycle of part Q is $T_q$, the failure warning time of the entire car is $T_0 = k*T_q - T_{min}$. When the car owner continues to use the car for $T_0$, the car control system reminds the owner that the car has been in the failure warning state in real time.

## 1.3 Implementation Program

1.3.1 System Architecture

According to Fig.1 and the introduction of the program implementation steps, it can be seen that the automotive fault diagnosis and early warning system consists of three modules: real-time information data acquisition, information data processing and expert system, and the functions of each module are as follows.

(1) Information Data Acquisition Layer

The premise of an automobile fault diagnosis and early warning system is to acquire the required accurate information in real time. The real-time information data acquisition module selects the corresponding sensors and installs them in suitable positions according to the system requirements, senses the relevant information data according to the requirements, and then applies IoT technology to transmit the measured information data to the information data processing module.

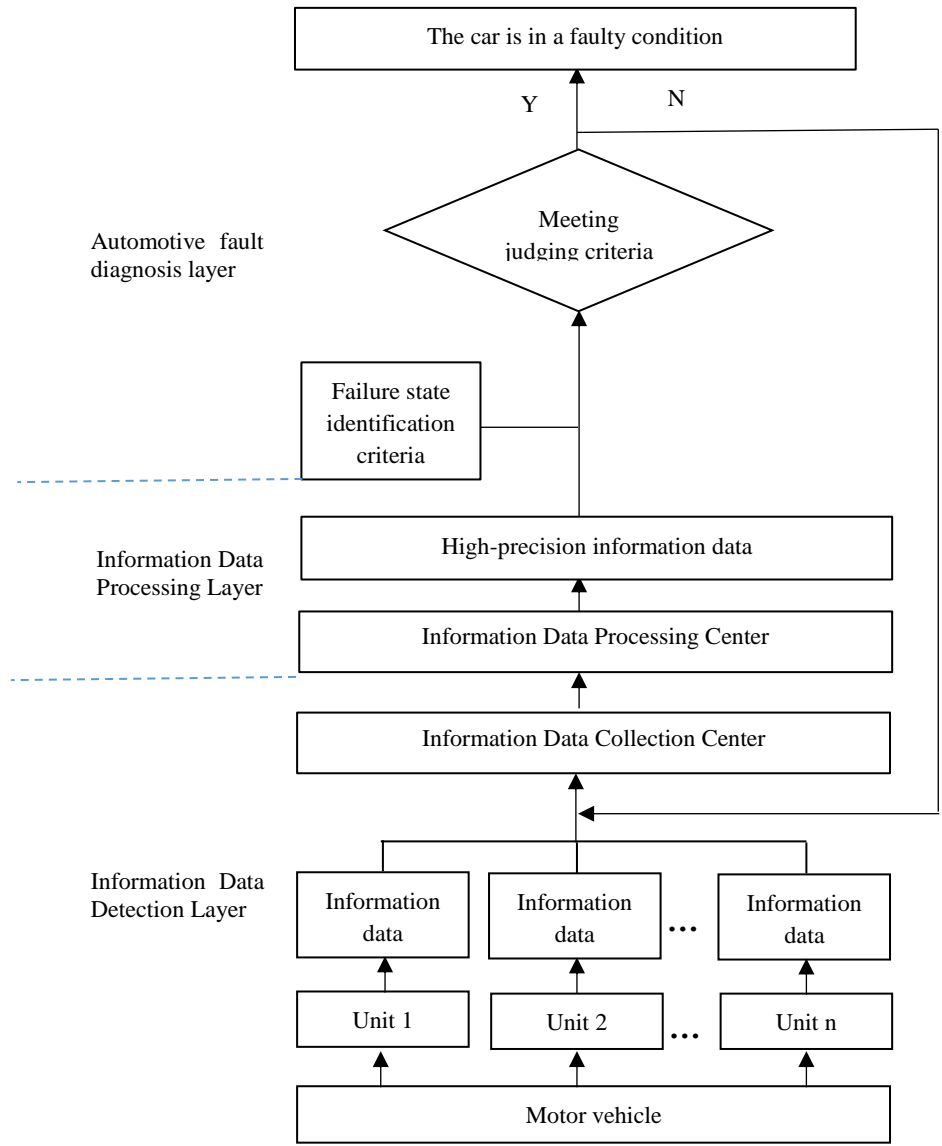

**Fig. 1** Workflow diagram of automobile fault diagnosis and warning system

(2) Information Data Processing Layer:

Owing to the variety and quantity of the collected information, the information data processing center. Therefore, it is necessary to categorize and store collected data for subsequent processing. Subject to the joint action of various factors, the collected information data cannot avoid unpredictable measurement errors. The fusion of various types of information data in the database improves the progress of information data detection and provides accurate information data for automotive fault diagnosis systems.

(3) Automobile fault diagnosis layer

Set the judgment standard for various types of automobile work unit failures and judge the current operation status of each work unit by comparing the acquired high-precision information data with the judgment standard for failure warning. If the measured value of the acquired information is equal to the judgment standard set by the system, it means that the car is already in the failure warning state; otherwise, it will continue to collect all types of information data generated by the work units, and the cyclic process of processing the information data and judging the failures until the car is in the failure warning state.

1.3.2 Technical difficulties faced in the implementation process

From the above, it can be seen that three problems are faced during the development and implementation of automotive fault diagnosis and early warning methods:

(1) The real-time nature of information data acquisition

The real-time nature of information collection is an important basis for measuring the efficiency of an information-detection system. When an automobile fault diagnosis and early warning system works, multiple types of massive information data must be collected, and there are large differences between different information data collection and transmission methods. However, with the development of electronic information and network technology, the technology to realize real-time collection of various types of information data has matured.

(2) Accuracy of information data value

Information data acquisition and data accuracy are affected by multiple factors, which are the main factors affecting the performance of the equipment, working environment, and signal transmission. Information data in the acquisition process, all kinds of influencing factors with changes in time and space, so always in dynamic change, the information data obtained by the detection system is bound to have unpredictable measurement errors. Therefore, eliminating these irregular measurement errors and improving the accuracy of information data have become technical challenges in the field of information data detection.

(3) The scientific nature of expert system decision-making

The expert system determines the automobile fault diagnosis and early warning system reliability of the work of the top building according to the information data provided by the information data detection system, accurate judgment of the car's current operation of the state, and diagnosis of the conclusion and proposed response measures. Therefore, the scientific nature of expert system decision making is an important indicator of the advanced reliability of the working system. Automotive fault diagnosis and early warning systems involving the function mode are not significant, and with current computer technology, the development of automotive fault diagnosis and early warning functions can be realized, and expert systems have no technical problems.

From the working principle and process of the engine fault diagnosis and early warning system, the accuracy of the information data is the entire system reliability of the premise, only to obtain high-quality information data to accurately determine the car's current operating state, accurate judgment of the type of failure, and the use of scientific preventive measures. However, because of the complexity of the working conditions of an automobile, the environment is variable, and the equipment operation of various types of information is inevitably subject to the common influence of multiple factors, bringing the irregularity of the measurement error, greatly affecting the automobile fault diagnosis and early warning system of the reliability of the work. Therefore, the development of detection technology for high-precision information data has become a difficult problem in the implementation of automotive fault diagnosis and early warning programs.

## 2. High-precision detection method of information data

Since the quality of information data affects the reliability of the work system and the scientificity of decision-making, it has been widely considered by all walks of life and has become a hot spot and a key area of research for most experts at present. Over the decades, with the joint efforts of many experts, promising research results have been obtained. According to existing research, measures to improve the quality of information data can be divided into two categories: hardware and software methods..

**2.1 Hardware method**

The essence of the hardware method is to improve the detection accuracy of the information data using high-performance detection instruments. The main research results for the hardware methods are as follows: Hu et al. [17] proposed a high-precision safety valve test architecture with three testing channels and effectively solved the problems of current safety valve testing. Huijun et al. [18] developed a comprehensive sliding-separation test platform for RV reducers to realize high-precision and high-display test performance for various RV reducer parameters. Reference [19-20] proposed a load differential radiation pulse on a transient electromagnetic high-performance radiation source for pulse-scanning detection to solve the problems of urban electromagnetic interference and insufficient harmonic components emitted by radiation sources. Reference [21-23] designed a hardware system based on radar and realized a real-time detection function for underground space-related information by enlarging detection information. Jiaqi et al. [24] proposed a one-stage remote sensing image object detection model: a multi-feature information complementary detector (MFICDet), which can improve the ability of the model to recognize long-distance dependent information and establish spatial–location relationships between features.

However, in engineering applications, we found that the hardware method had the following shortcomings.

① The hardware method to improve the accuracy of information data is to improve the detection performance of the detection equipment, but the detection accuracy of the detection equipment and its price are directly proportional to the serious shortcomings of the low cost.

② The detection equipment must improve the strength and reduce the signal detection error performance, but if you want to have this function at the same time, you need to use two kinds of equipment, which not only increases the cost but also increases the use of space, but is also prone to performance mismatches due to system instability problems.

**2.2 Software method**

In recent years, most researchers have attempted to use software methods to achieve high-precision information-data detection to solve the shortage of hardware methods for information-data detection in complex environments. The essence of the software method is the information data-fusion algorithm. To date, many research results have been obtained, such as fuzzy set theory [25], fuzzy neural networks [26], probability models [27], and particle swarm optimization algorithms [28]. For example, Huo et al. [29] proposed an integral infinite log-ratio algorithm (IILRA) and an integral infinity log-ratio algorithm based on signal-to-noise ratio (BSNR-IILRA) to improve the detection accuracy of the laser communication detection position in the atmosphere. Zhiyuan et al. [30] proposed a normalized-variance detection method based on compression sensing measurements of received signals and solved the problem of fast and accurate spectrum sensing technology under the condition of a low signal-to-noise ratio. Liu et al. [31] proposed a target detection algorithm based on the improved RetinaNet, which is suitable for transmission-line defect detection and improves the intelligent detection accuracy of UAV in power systems. Cheng et al. [32] proposed a lightweight ECA-YOLOX-Tiny model by embedding an efficient channel attention (ECA) module, which has a higher response rate for decision areas and special backgrounds, such as overlapping small target insulators, insulators obscured by tower poles, or insulators with high-similarity backgrounds. Liu Wenqiang et al. [33] introduced a point cloud segmentation and recognition method based on three-dimensional convolutional neural networks (3-D CNNs) to determine the different components of the catenary cantilever devices. Yin et al. [34] proposed a complementary symmetric geometry-free (CSGF) method that makes the detection of cycle slips more comprehensive and accurate. Lingfeng et al. [35] established a junction temperature model based on a multiple linear stepwise regression algorithm, and used it to extract high-precision intersection online temperatures. However, through the analysis of various current software methods, the following deficiencies were found in the detection of information in complex environments:

① The application of these algorithms is carried out in the ideal environment of the laboratory, failing to take into account that the detection of information in engineering practice will be affected by various types of dynamic changes in the boundary factors, thus bringing the information data measurement value presents the reality of irregular measurement errors, so there is a feasibility and reliability of engineering applications.

② These algorithms have the function of improving the detection accuracy of information data, but fail to take into account the information data in the transmission process of signal interference and energy loss. Signal interference can easily lead to the detection of signal distortion phenomenon in the transmission process, and the energy loss of the measured value of information data is low.

Therefore, they have a limited application area; therefore, an ideal high-precision detection method for information under the joint action of multiple influencing factors has not been developed. To solve these problems, our team has been using the method of fractional calculus theory in data processing for many years [16,36-38] and has found that fractional differential operators are suitable for studying non-linear, non-causal, and non-stationary signals and have dual functions of improving detection information and enhancing signal strength. Therefore, this paper proposes the application of the fractional-order differential data fusion algorithm to explore the measurement of high-quality information data in complex environments under established quality standards.

**3. information data fusion algorithm based on fractional order differential operator**

**3.1 Fractional-order Calculus Definition**

Fractional calculus is a branch of calculus derived from the theory of integer order calculus, which first appeared in L'H δ- in the letter written by pit to Leibniz [39]. The theory of fractional calculus was established 300 years ago, but it has long been in the stage of purely mathematical theoretical analysis and derivation by mathematicians. In recent years, renowned mathematicians such as Liouville have begun to focus on the improvement and development of fractional calculus theory and have established a basic system architecture for fractional calculus theory. The fractional derivative is actually any order fractional calculus, which is an important branch of mathematics that has just been developed from the n-th derivative and n-th integral in recent years. Although it has a history of over 300 years, owing to the different research objects in various application fields, fractional calculus has always lacked a unified definition in various fields. Currently, there are three well-known definitions of fractional derivatives in basic theory and engineering application research: Grunwald Letnikov, Riemann Liouville, and Caputo Riesz [40-41].

3.1.1 Grunwald-Letnikov Fractional Order Differential Definition

For any real number v, the integer part of $v$ is denoted by $[v]$. Assuming that the function has $n + 1$ continuous derivatives on the interval $[v, t]$, and when $v > 0$, $n \geq [v]$, the fractional-order $v$ derivative is defined as

$$_aD_t^v f(t) = \lim_{h \to 0} f_h^{(v)}(t) = \lim_{h \to 0} h^{-v} \sum_{j=0}^{\left[\frac{t-v}{h}\right]} (-1)^j \binom{v}{j} f(t - jh) \qquad (2)$$

Where $0 \leq n-1 < v < n$; $\binom{v}{j}$ is the Binomial coefficient, and

$$\binom{v}{j} = \frac{v(v-1)(v-2)....(v-j+1)}{j!} \qquad (3)$$

This definition is widely used in the field of numerical calculation and has evolved from the classical definition of the integral derivative of a continuous function. If the derivative order in (1) is extended to a fraction, the definition of the Grunwald Letnikov fractional calculus can be obtained.

3.1.2 Riemann-Liouville Fractional Calculus Definition

According to Cauchy Antiderivative formula, we can get:

$$_aD_t^v f(t) = \frac{d^n}{dt^n}\left[_aD^{v-n}f(t)\right] = \frac{1}{\Gamma(n-v)} \frac{d^n}{dt^n}\left[\int_a^t (t-\tau)^{n-v-1} f(\tau) d\tau\right] \qquad (4)$$

Where $0 \leq n-1 < v < n$

This definition is an improved expression of Grunwald-Letnikov's definition, which simplifies the calculation process of fractional calculus based on its basic properties of fractional calculus. Riemann–Liouville fractional calculus can be obtained from (3).

3.1.3 Definition of Caputo Fractional Calculus

In fact, the Caputo fractional calculus definition is the same as the R-L definition. The equation is:

$$_a^C D_t^v f(t) = \frac{1}{\Gamma(n-v)} \int_v^t \frac{f^n(\tau) d\tau}{(t-\tau)^{v+1-n}} \qquad (5)$$

Where $0 \leq n-1 < v < n$.

With regard to the three classical definitions of fractional order calculus, many scholars at home and abroad have studied and analyzed a variety of practical problems from different application perspectives and have obtained different forms of expressions for the definitions of fractional order calculus. The Riemann-Liouville definition of fractional order calculus and Caputo definition of fractional order calculus are both improvements of the Gnimwald-Letnikov definition of fractional order calculus. The Grumwald-Letnikov definition of fractional-order calculus can be converted to convolution operations in numerical realizations, which makes it very suitable for signal processing applications. The Riemann-Liouville definition of fractional-order calculus is mainly used to compute the analytic solutions of simpler functions, whereas the Caputo definition of fractional-order calculus is suitable for problems of initial marginal values of differential equations of fractional order. The Caputo definition of fractional-order calculus is suitable for analyzing the first marginal value problems of fractional-order differential equations, and is therefore well suited for engineering applications. The Caputo definition of fractional-order calculus is suitable for analyzing the initial margin problems of fractional-order differential equations and is therefore suitable for engineering applications.

From the above discussion, it can be seen that there are certain connections and differences between the three definitions of fractional calculus. After analyzing the three definitions of fractional calculus, it can be seen that compared with the other two definitions, the G-L definition is widely used by the engineering community because of its simple and efficient calculation process, although there is a slight lack of computational accuracy. Therefore, this study applies the fractional-order differential algorithm under the definition of G-L to study high-precision detection technology for various types of automobile information data.

**3.2 Application of fractional order differentiation to signal processing.**

Through the study of the spectral characteristics of fractional-order differential operators, many experts have found that fractional-order differential operators are suitable for the study of non-linear, non-causal, non-smooth, and other characteristics of uncertain signals and have been applied to electrochemistry, mechanics of materials, fluid field theory, electromagnetic field theory, biomedicine, signal processing, and other subject areas. In engineering practice, because the information data in the detection process is affected by the performance of the equipment and the impact of the working environment, resulting in the collection of information data, there are no rules to follow the non-linear changes in the state of the signal transmission process that will be subject to various types of signal interference and energy loss under the influence of unpredictable signal distortion problems. Therefore, most of the signals obtained in engineering practice are non-linear, noncausal, and non-smooth. This study intends to apply fractional-order differential operators to the processing of engineering detection signals based on the study of the spectral characteristics of fractional-order differential operators.

Assuming that a square-integrable energy signal $S(t)$, and $S(t) \in L^2(R)$, $S(t)$ is captured by applying Internet of Things (IoT) technology in engineering practice, a Fourier transform can be applied to the signal to obtain.

$$S(\hat{\omega}) = \int_R S(t) e^{-i\omega t} dt \quad (6)$$

According to the nature of the Fourier transform, the v-order differential operator is equal to the v-order differential multiplicative function $\hat{d}(\omega) = (i\omega)^v$ of the multiplicative operator, the fractional-order derivatives as $S^v(t)$, and the following equation can be obtained:

Make,

$$d^v(\hat{\omega}) = |\omega|^v e^{i\theta^v(\omega)}, \theta^v(\hat{\omega}) = \frac{v\pi}{2} \operatorname{sgn}(\omega) \quad (7)$$

Then the (6) can be simplified as:

$$D^v S(t) \overset{FT}{\Leftrightarrow} (\hat{D}S)^v(\omega) = (i\omega)^v S(\hat{\omega}) = d^v(\hat{\omega}) S(\hat{\omega}) \quad (8)$$

Therefore, from the perspective of signal modulation, the physical meaning of detecting the fractional differentiation of signals is equivalent to generalized amplitude and phase modulation. From the perspective of signal processing, the v-order fractional calculus operation of the detected signal is equivalent to establishing a linear time-invariant filtering system for the signal. The filtering function of the time-invariant filtering system is $d^v(\hat{\omega}) = (i\omega)^v = |\omega|^v e^{i\theta_v(\omega)}$.

Fig. 1 shows the amplitude–frequency characteristics of the fractional differential operator when the fractional order $v$ takes different values. From Fig.1, it can be observed that

① In the processing of detection signals, the fractional-order differentiation operation can enhance the middle- and high-frequency parts of the signal, and the enhancement amplitude grows sharply and non-linearly with the increase in the frequency and fractional-order differentiation order. Therefore, fractional-order differentiation enhances the high-frequency components of the signal while non-linearly preserving the very low-frequency components of the signal.

② With the increase in fractional order and signal frequency, the fractional-order differential operator significantly enhances the middle- and high-frequency parts of the energy signal. In the high-frequency region of the signal, the data difference between the signal enhancement values of the fractional-order differential operators at different orders decreases, and with an increase in frequency, the differential operators at the same order have essentially the same effect on the enhancement of the signal strength at different frequencies.

According to the relationship between signal strength and the fractional differential operator shown in Fig.1, scientists found that fractional calculus under different differential operators has different enhancement effects on signals; for example, lower frequency components can be kept non-linear. From a physical perspective, the fractional differential processing of signals can be understood as a generalized amplitude phase modulation. Because the fractional differential operator has the function of significantly strengthening its signal strength when processing high-frequency signals, it can significantly improve the strength of high-frequency and low-frequency signals. Therefore, in the fractional-order differential operator signal denoising process, not only can the key information of the signal edge area be maintained, but the basic information of the signal middle area can also be well preserved, which can significantly improve the signal strength.

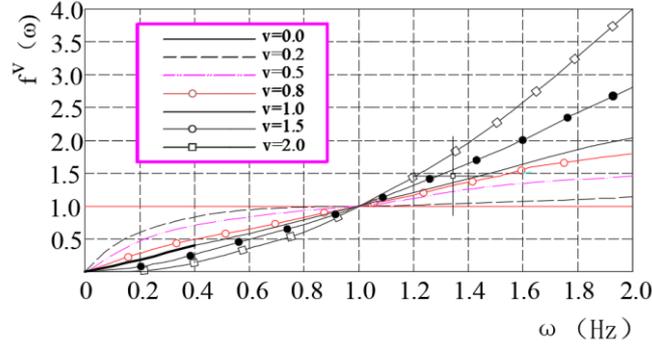

**Fiq.2** The Amplitude-frequency Curve of the fractional differential operator

## 4. Fractional order differential operator based automobile information data fusion processing technology

### 4.1 Information data fusion processing model.

In the expression of the fractional-order calculus under the G-L definition, the frequency ω is the influence factor of the signal $f(\omega)$, Signal $f(\omega)$, and the differential result is $f^v(\omega)$, where $\omega \in [\omega_1, \omega_2]$, we can adjust the frequency ω. If the measured value of is regarded as the corresponding variable $x$, then $\omega_i = x_i$, and, where $x \in [a, b]$, the function $f(x)$ is a function of the influence factor x. Referring to Fig.2, we can obtain a schematic of the fractional-order differential function $f^v(x)$ with the variable function $f(x)$, as shown in Fig.3.

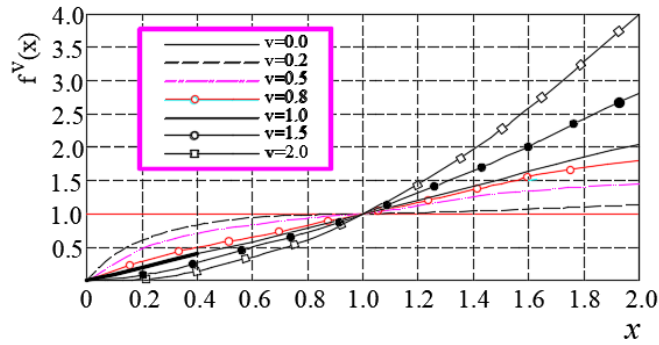

**Fig. 3** Schematic diagram of the effect of the function $f(x)$ processed by fractional order differentiation

From Fig.3, we can see that after the function $f(x)$ is processed by fractional-order differentiation, the relationship between the value $f^v(x)$ and the influencing factor $x$ is as follows:

① The fractional-order differential operator plays an auxiliary role in the signal data processing. After the signal data are processed by the fractional-order differential operator, the value of the information data increases as the values of the fractional order v and influence factor x increase.

② When $v=0$, the fractional differential operator does not amplify the measured value of the information data.

③ When $0 < v < 1$, and at the same differential order v, the fused information data values are amplified as the value of the influence factor x increases, but the growth of the amplification factor gradually decreases.

④ When $0 < v < 1$, at the same value of the influence factor x, the amplification coefficient of the fractional-order differential operator on the information data value increases with an increase in the fractional order $v$, but the gap between the amplification coefficients of the different orders on the information data gradually decreases with an increase in the influence factor $x$.

According to the Grunwald-Letnikov definition, the duration interval $(a,b)$ of the influence factor $x$ of the functional equation $E(x)$ of the information data collected by the IoT is homogenized in steps $h$, which yields. This yields the fractional-order differential equation of the function $E(x)$ with respect to the influence factor $x$ under the Grunwald-Letnikov definition:

$$_a^G D_t^v E(x) = \lim_{h \to 0} E_h^v(x) = h^{-v} \sum_{i=0}^{\left[\frac{b-a}{h}\right]} (-1)^i \binom{v}{i} E(x - ih) \qquad (9)$$

Where $0 \leq n - 1 < n$, $\binom{v}{i}$ are binomial coefficients and

$$\binom{v}{i} = \frac{v(v-1)(v-2)\ldots[v-(i-1)]}{i!} = \begin{bmatrix} -v \\ i \end{bmatrix} \qquad (10)$$

Substituting (10) into (9) yields the model of the information data fusion algorithm based on fractional order differentiation under the Grunwald–Letnikov definition:

$$_{a}^{G}D_{t}^{v}E(x) = \lim_{h \to 0} E_{h}^{v}(x) = h^{-v}\sum_{i=0}^{\left[\frac{b-a}{h}\right]}(-1)^{i}\begin{bmatrix} -v \\ i \end{bmatrix}E(x-ih) = h^{-v}\sum_{i=0}^{\left[\frac{b-a}{h}\right]}w_{i}^{v}E(x-ih) \qquad (11)$$

where $w_i^v$ is the value of the weighting coefficient, calculated using the formula

$$w_0^v = 1, w_i^v = \left[\frac{i-(v-1)}{i}\right]w_{i-1}^v, i = 1, 2\cdots \qquad (12)$$

**4.2 Detection information data fusion processing**

Adopting the integral idea under the distributed system, the production information monitoring system based on Internet of Things (IoT) technology first categorizes the collected information and then judges the validity of similar information data. If the deviation value between similar data collected exceeds the tolerance range set by the system, the multi-sensor detection data fusion model based on the fractional-order differential operator is applied to fuse the data, and efficient fusion between similar data is finally realized. Therefore, the fusion process of monitoring information data under the IoT based on the fractional-order differential operator is as follows:

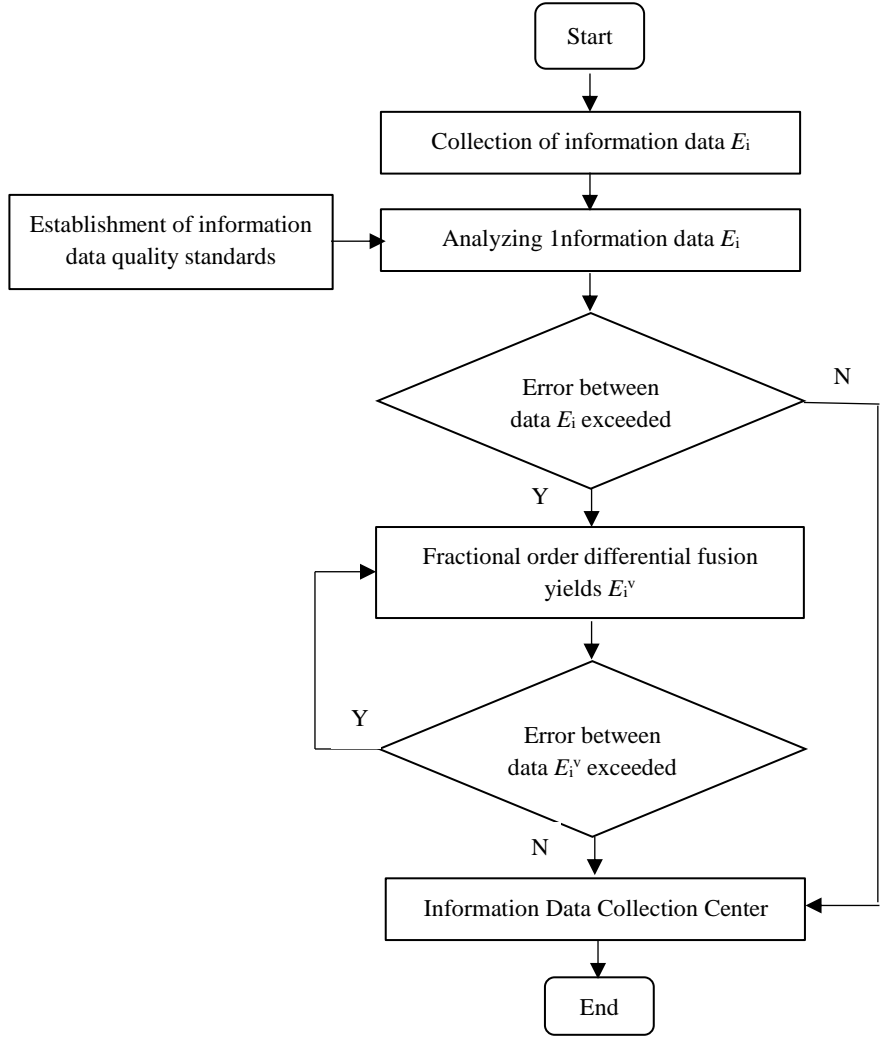

**Fig.4** Detection fusion procession

Step 1: IoT technology is applied to collect the required information data, which is affected by factors such as equipment performance, working environment, and signal transmission, and there is an unpredictable measurement error between the information data Ei collected by the information detection system.

Step 2: Setting the information collection quality judging standard, judging whether the error between the collected information data is within the error range set by the detection system; if it is within the error range, the detection results are sent directly to the data collection center; if the error between the data exceeds the error range set by the system, then proceeding to the next step.

Step 3: According to the data attributes of the detected information, select an influence factor x of the measured value of the detected information, synthesize the measured value $E_i$ and the corresponding influence factor $x_i$ of each sensor, and apply MATLAB software to fit a functional relationship equation $E(x)$ between the measured value $E_i$ and the influence factor $x_i$.

Step 4: With reference to the detection data fusion model (11), construct a fractional-order differential data fusion model of equation $E(x)$, and apply MATLAB software to calculate the fractional-order differential processing results $E^v_i$ of each detection data $E_i$.

Step 5: Analyze and compare the error between the fractional-order differential processing results of each detection data Ei, and evaluate the fusion results. If the error between the results is within the range set by the system, fused data $E^v_i$ is sent to the data center. Otherwise, return to step 3 as the initial value and apply the fractional-order differential operator fusion processing again until the error between the final fusion results is within the range set by the system.

Step 6: The fractional-order differential fusion processing of the detection data of the detection information ends.

## 5. Application examples

### 5.1 Experimental platform

Automobile structures are complex; therefore, there are many types of faults, and it is difficult to identify their characteristics. The common faults of each part of the car include common faults of electronically controlled gasoline injection engines, common faults of the transmission system, common faults of the steering system, common faults of the braking system, common faults of the walking system, and electrical system. Owing to space constraints, this study can only select the important devices in the car engine and body as the object of study to explore the feasibility and effectiveness of the method described in the paper applied in engineering practice.

From the perspective of diagnostic accuracy, a car close to the fault state was selected as the research object. After careful selection, this study takes a Buick Insight 2013 GT1.6L model car as the experimental object, the vehicle service life is 9 years and 7 months, and the use of the engine and the body have undergone minor repairs during the period of use. During the experiment, the car operates normally; at the same point in time, the corresponding testing equipment is applied to collect the main parameter data of typical faults, analyze the data to judge the operating status of the engine and gearbox, predict the time of fault occurrence, synthesize the results of the fault diagnosis judgment of the engine and body, and analyze the current operating status of the entire car and the time of fault occurrence.

### 5.2 Operational status evaluation technology for automobile engines

For an electronically controlled gasoline injection engine, its important faults include inability to start or start, poor idling, rattling, high fuel consumption, high lubricating oil consumption, abnormal smoke from the exhaust pipe, weak acceleration, oil seepage, and water leakage. The fault phenomena are mainly rattling, oil temperature abnormalities, and water temperature abnormalities. Among them, rattling is the main expression of engine failure, and an investigation found that approximately 70% of automobile engine failures are manifested through rattling. Therefore, this study chose rattling as a measure of engine failure to study engine Real-time operational status evaluation technology.

With the growth in the use of automobile engines, improper operation, maintenance quality, and the influence of the natural environment, various parts due to wear and tear, breakage, loosening, aging, poor contact, short-circuit, and disconnection, are produced more than the prescribed rattles in the work, such as knocking sound, overspeed running whistling, parts rubbing sound, and so on. Therefore, rattling is an unavoidable fault phenomenon in the engine life cycle. Fault diagnosis can be used to determine the general rules and characteristics of a fault from the most intuitive manifestation of rattling, which will bring great convenience to the diagnosis of automobile engine faults. Therefore, obtaining accurate engine rattling measurement data and determining the change rule of rattling can accurately diagnose the type of engine failure and operation status and predict the time of failure.

5.2.1 Information Data Detection

To simulate a real working environment, the experiment was conducted in an outdoor environment, which was more than 5 m away from the surrounding buildings, and the height was approximately 1.5 meters. To facilitate the comparison of the fractional-order differential algorithm and other instruments between the information data detection effect, the experiment used nine IFM capacitive sensors to detect the noise value of the engine at different locations: the sensors were distributed in the engine and the left and right sides, each side of the three, and were used to distribute the distance from the engine is 0.3 meters. When the engine was running for 10 min after the detection of engine noise measurements, the experiment was sampled five times at 30 s intervals. The measured values are listed in Table 1.

Table 1. Measured values of vehicle engine noise information data (dB)

| Sensor number | No. of measurements | | | | | Mean value $F_i$ | standard deviation $S_i$ |
|---|---|---|---|---|---|---|---|
| | NO.1 | NO.2 | NO.3 | NO.4 | NO.5 | | |
| F1# | 70.4 | 70.6 | 71.4 | 78.8 | 70.9 | 70.83 | 0.38 |
| F2# | 71.1 | 71.3 | 71.2 | 71.5 | 71.4 | 71.30 | 0.14 |
| F3# | 70.9 | 71.2 | 70.8 | 70.6 | 79.28 | 70.88 | 0.22 |
| F4# | 70.7 | 71.5 | 70.9 | 71.2 | 70.9 | 71.04 | 0.28 |
| F5# | 71.1 | 70.8 | 70.6 | 71.4 | 71.5 | 71.29 | 0.34 |
| F6# | 70.8 | 71.4 | 71.2 | 71.2 | 70.8 | 71.08 | 0.24 |
| F7# | 70.7 | 71.1 | 62.1 | 71.6 | 71.1 | 71.13 | 0.32 |
| F8# | 71.2 | 71.4 | 71.4 | 71.1 | 70.6 | 71.14 | 0.29 |
| F9# | 70.8 | 71.6 | 70.7 | 70.5 | 70.8 | 70.88 | 0.38 |
| Average $F$ | | | | | 71.038 | | |
| Standard deviation $S$ | | | | | 0.144 | | |

5.2.2 Information Data Analysis

From the experimental environment above, it can be seen that the measured value of the engine rattle received by the data center in the sensor acquisition will be affected by the sensor performance, working environment, and other factors; in the process of information transmission, it will be subjected to complex signal interference and energy attenuation caused by the detection error. Because the nine sensors were placed in the same place, the working environment of each sensor and the signal interference during data transmission were essentially the same. Therefore, the variability between the detected values of each sensor is mainly related to the performance of the detection instrument, and the error between the measurements is mainly related to the changes in the working environment at different time points. However, because the time difference between the five measurements in this experiment was relatively short, and the possibility of a sudden change in the working environment within 2.5 minutes was relatively small, it can be assumed that the variability between the measured values of the different sensors in this experiment was mainly due to the variability in the performance of the detection instrument.

In the experiment, we set the measurement error of the experimental data to 10%, that is, the qualifying criterion of the data was $64 < F < 78$. Therefore, the data 78.8, 79.2 and 62.1 in Table 1, are obviously out of error and should be removed. This can be seen from the data shown in Table 1: the average value of five measurements is taken as the measurement result of each sensor, and the average value of all measurements is taken as the true value of the measurement $E$. A summary of the measured values of the engine noise data (as shown in Table 2) and the distribution of the detected data (as shown in Fig. 2) can be obtained. It can be observed from Table 2 that there is a large variability in the performance of the testing equipment, resulting in a large measurement error between the information data measured by each sensor. From the detection data distribution curve of each sensor shown in Fig.2, it can be seen that the detection data do not present the characteristics of random distribution of the measured value in the laboratory near the true value, but rather in the measurement of the true value near the state of irregular discrete distribution. The above two points can be observed: affected by the performance of the testing equipment, the detection value of the engine rattle is seriously distorted, and it is difficult to provide accurate detection values for fault diagnosis and early warning systems, which seriously affects the reliability of the system.

Table 2. Summary of experimental data on vehicle engine rattles (dB)

| Sensor number | F₁# | F₂# | F₃# | F₄# | F₅# | F₆# | F₇# | F₈# | F₉# |
|---|---|---|---|---|---|---|---|---|---|
| Average value $F_i$ | 70.83 | 71.30 | 70.88 | 71.04 | 71.08 | 71.08 | 71.13 | 71.14 | 70.88 |
| Standard deviation $S_i$ | 0.38 | 0.14 | 0.22 | 0.28 | 0.34 | 0.24 | 0.32 | 0.29 | 0.38 |
| Measured true value $F$ | | | | | 71.038 | | | | |

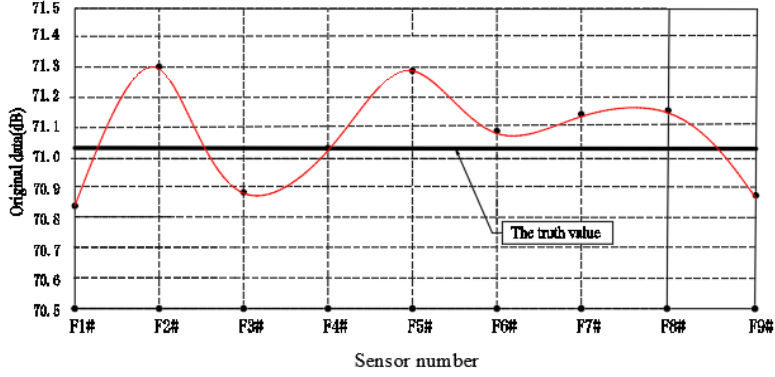

**Fig.5** Detection data distribution curve

### 5.2.3 Information data fusion processing
(1) Selection of the detection value influence factor

As mentioned above, nine sensors work in the same working environment, at the same time, and at the same frequency to detect the value of the engine rattle; thus, the work of the environmental impact and signal interference in the transmission of information is basically the same, and the detection accuracy between the sensors is mainly related to the performance of the detection equipment. In the measurement technology, the standard deviation $S$ is the best parameter to measure the performance of the detection equipment; therefore, this case adopts the square standard deviation S as the influencing factor of the sensor's measurement value, and analyzes the correlation between the performance of the detection equipment and the measurement value. The measured values of each sensor and their measurement standard deviations ($S_i$) are listed in Table 2.

(2) Calculation of the relationship between the measured values and influencing factors

Because the least-squares method has the advantage of not requiring a priori data in the data processing process, it is widely used in fitting the polynomials of the function of the detection data and can achieve ideal data fusion accuracy; therefore, the least squares method is used in this case to fit the function $F(x)$ between the measured value $F_i$ of each sensor and its influence factor $x_i$. The parameter values of the measured value $F_i$ and its influence factor $x_i$ (standard deviation $S$) for each sensor are listed in Table 2. To find the functional equation $F(x)$ between the measured value $F_i$ and its influence factor $x_i$, the mathematical expression of the function $F(x)$ is assumed:

$$F(x) = a_0 + a_1 x + a_2 x^2 + \cdots + a_n x^n \quad (13)$$

From the above equation, it can be seen that the essence of obtaining the function $F(x)$ is the calculation of polynomial order $n$ and coefficient ai.

According to the definition of the method of least squares, the polynomial order of the fit should be less than the number of data samples; therefore, the order $n$ of the function $F(x)$ sought in the case is less than 9. Because the fitting accuracy of the function $F(x)$ does not improve with an increase in the order $n$, the optimal order value must be found between $0 < n < 9$ so that the error between the function $F(x)$ and the true value of the measurement is minimized. By applying the Polyfit function in MATLAB software, the total measurement error value of (11) under different orders can be obtained, and by comparing the total measurement error under different orders, it can be seen that when the fitting equation is quadratic, the total error between the fitting value of each sensor and the true value of the measurement is minimized, with a value of 0.1153. At this time, the equation of the $F(x)$ function is

$$F(x) = a_0 + a_1 x + a_2 x^2 \quad (14)$$

By applying the Polyfit function in MATLAB math software, it can be concluded that $a_0$=71.37465, $a_1$=-0.9140, and $a_2$=-0.4923, so that the functional equation between the detection data $F_i$ and the influence factor $x_i$ based on the least squares method can be concluded as follows:

$$F(x) = -0.4923 * x^2 - 0.9140 * x + 71.347 \quad (15)$$

(3) Detection data fusion model based on fractional-order differential operator

The data fusion algorithm model of coal mine detection information based on the fractional-order differential operator under Grunwald-Letnikov's definition is shown in (9). If we want to obtain a detailed mathematical model, we should first determine the values of order $v$ and step size $h$ in (9).

① Selection of order $v$

From the characteristics of the fractional-order differential operator at different orders v shown in Fig. 1, it can be seen that in the low-frequency part, the fractional-order differential operator has a certain enhancement effect compared to the integer-order operator when the differential order $0 < v < 1$, and the enhancement effect becomes more obvious with a decrease in order, which has the effect of non-linearly preserving the low-frequency components of the signal. To save space, the intermediate value $v=0.5$ of [0,1] is used in this case to explore the fusion processing effect of the fractional-order differential operator on the detection data.

② Value of step size h

According to the Grunwald-Letnikov definition of fractional-order calculus, the smaller the value of the step size $h$ is, the more accurate the calculation results; however, it will also bring the problem of increase the calculation workload. Therefore, the value of the step size $h$ should be considered for both calculation accuracy and speed. In this case, the value range of the detection data influence factor $x \in [0.14, 0.38]$; therefore, the step size $h$ can be taken as 0.01, so that the number of steps to be calculated is $i = 24$, which can simultaneously meet the requirements of both calculation accuracy and speed.

Now, the fractional order $v = 0.5$ and the step size $h = 0.01$ are substituted into (9), and the fusion algorithm model of the detection data in this case can be obtained:

$$F^{v}(x) \approx \lim_{h \to 0} h^{-v} \sum_{i=0}^{[b-a]/h} \begin{bmatrix} -v \\ i \end{bmatrix} F(x-ih) = 0.01^{-0.5} \left[ F(x) - 0.5F(x-0.01) - \frac{0.5(1-0.5)}{2!} E(x-2\times0.5) - ... \right.$$
$$\left. - \frac{0.5(1-0.5)(2-0.5)..(24-1-0.5)}{24!} F(x-0.24) \right] \quad (16)$$

5.2.4 Information data fusion results

Combining (13) can derive the value of $F(x-ih)$ at each step and substituting the calculation results into (16), the 0.5-order differential fusion processing results of each sensor's measured value $F(x_i)$ with respect to its influence factor $x_i$ can be obtained, as shown in Table 4.

Table 3  0.5-order differential fusion values of measured data from each sensor

| Sensor number | F1# | F2# | F3# | F4# | F5# | F6# | F7# | F8# | F9# |
|---|---|---|---|---|---|---|---|---|---|
| Post-fusion value $F_i^{0.5}$ | 112.2 | 112.31 | 112.28 | 112.26 | 112.27 | 112.15 | 112.31 | 112.29 | 112.2 |
| Standard deviation S | 0.38 | 0.14 | 0.22 | 0.28 | 0.34 | 0.24 | 0.32 | 0.29 | 0.38 |
| Pre-fusion value $F$ | 71.038 | | | | | | | | |
| Enlargement ratio | K=112.252/71.038=1.582 | | | | | | | | |

Table 4  0.5-order differential fusion values of measured data from each sensor

| Sensor number | F1# | F2# | F3# | F4# | F5# | F6# | F7# | F8# | F9# |
|---|---|---|---|---|---|---|---|---|---|
| Average value $\overline{F_i^{0.5}}$ | 71.01 | 71.11 | 71.06 | 71.05 | 71.06 | 70.98 | 71.08 | 71.07 | 71.01 |
| Standard deviation after | 0.034 | | | | | | | | |

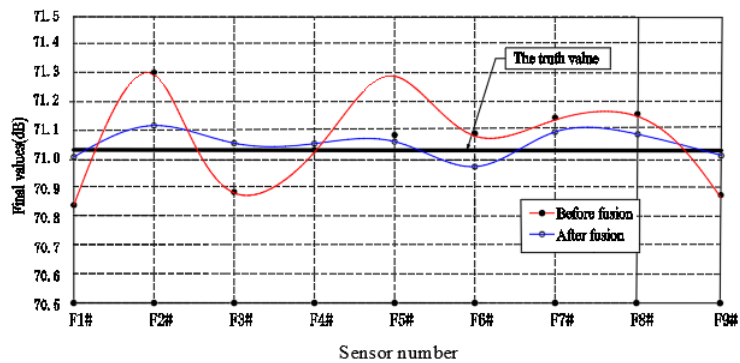

Fig.6 Comparison of distribution curves before and after information data fusion

From the final fusion results shown in Table 4, it can be seen that after 0.5-order fractional-order differential operator fusion processing, the detection value of the engine rattles was improved, and the average value of the measurements after fusion processing was 112.252, which not only did not decrease, but also improved the intensity by 58.2% compared with the true value of the measurements. After fusion, the error between the detection values of each sensor is smaller, and the dispersion between the original data decreases dramatically, as can be seen from the distribution curve of the fused detection data shown in Fig. 6: After 0.5 order fractional order differential operator fusion processing, the engine rattling measurement values are randomly distributed near the average value, and the detection accuracy is improved by more than 5 times, and the measurement accuracy S between the fused information data is substantially lower than 0.075 of the algorithm used in the literature [29], 0.043 of the algorithm used in the literature [30] and 0.054 of the algorithm used in the literature [31], and what's more important is that the fractional order differential algorithm can enhance the strength of the information data, enhance the information anti-interference ability in the measurement and compensate the measurement error caused by energy loss..

5.2.5 Operational status evaluation

The noise during the normal operation of the car mainly comes from the sound of the engine; therefore, this study considers the sound of the engine as the signature component of the car noise measurement value to judge the faulty operation status of the car.

(1) Fault diagnosis

In this experiment, the engine was selected from the normal operation of the car, in line with the motor vehicle noise measurement standards for sample car selection, in line with the requirements of motor vehicle noise test preparation. The experimental platform is located outdoors, and the surrounding buildings are more than 5 m in height and approximately 1.5 meters. The measurement of a total of nine sensors, and sensors in the engine, three sides of the uniform distribution, meet the motor vehicle noise test labeling requirements for the measurement area and sensor distribution. The number of information data acquisitions for each side of the engine during measurement is five, which meets the requirement of the standard that the number of measurements for each side is not less than. Therefore, the conditions of this experiment met the international requirements for motor vehicle noise-detection conditions.

According to China's motor vehicles driving outside noise extreme value $F_b \leq 74$ dB regulations, to ensure a safe period, this study will run the car warning coefficient $k$ set to 0.98, then $F_y = k * F_b = 72.52$. Referring to the accurate experimental data obtained in the experiment, $F_y = 71.038$, it can be determined that the current engine is in a normal working condition.

(2) Fault warnings

Because sound predicts an important parameter of engine failure, a failure warning is the acquisition of the point in time when the sound measurements are at the warning value set by the system. Therefore, the car measurement value $J$ is a function $F(t)$ at time $t$. The premise of accurate early warning of car faults is to accurately obtain the function $F(t)$ and then derive the corresponding warning time point $t_y$ according to the system-set early warning acoustic measurement value $F_y$.

In this study, we fitted the function $F(t)$ according to the measured value $F_i$ of the engine sound at different time points ti. It is reasonable that the more samples that are adopted, the more accurate the functional equation $F(t)$ is obtained, but this will lead to the shortcomings of high cost and low efficiency. Therefore, the number of samples used in the experiment was 6. The experimental samples were tested 0.5 months apart, and the experimental environment, information data collection, and processing methods were the same. The correspondence between the measured sample values $F_i$ and the time is presented in Table 5.

**Table 5.** 0.5-order differential fusion values of measured data from each sensor

| Sample value $F_i$（dB） | 71.038 | 71.086 | 71.173 | 71.335 | 71.624 | 72.207 |
|---|---|---|---|---|---|---|
| Point in time $t_i$ (months) | 0 | 0.5 | 1.0 | 1.5 | 2.0 | 2.5 |

By applying the Polyfit function in MATLAB software, we can obtain the total measurement error value of the eigenvalue under different orders, and by comparing the total measurement error under different orders, we know that when the fitting equation is a binary quadratic equation, the total error between the fitted value of the acoustic sound at each measurement point and the true value of the measurement is the smallest, with a value of 0.012. This can lead to detection based on the least-squares method. The functional relationship between data $F_i$ and influence factor $t_i$ is given by

$$F(t) = 0.2488*t^2 - 0.1865*t + 71.0735 \qquad (17)$$

Combining the acoustic measurements $J_y = 72.52$ at the time of the engine fault warning and (17), we obtain

$$0.2488*t^2 - 0.1865*t - 1.4465 = 0 \qquad (18)$$

Solving (18) shows that the engine failure warning time $t_y = 22.63$, i.e. about 22.63 months from the time of the first experiment.

**5.3 Real-time operational status evaluation technology of the vehicle body**

With the increase in vehicle use time, the vehicle in the process of traveling experiences a phenomenon of shaking, which is increasingly obvious. This type of shaking is a warning signal for vehicle problems. The causes are divided into four types: excessive engine carbon, ignition system problems, unstable oil pressure, and engine foot aging. Regardless of the cause, it is necessary to test this in time to eliminate serious potential problems.

The body of the Buick Insight was a load-bearing body. In this structure, the body not only serves as the mounting base for the engine and chassis assemblies but also functions as a frame, bearing all the loads. The operation is characterized by smooth driving, low inherent frequency, low noise, and low weight. Therefore, the main feature for determining the failure state of the body is the vibration frequency.

5.3.1 Information data acquisition

Owing to the diagnosis of a vehicle's operating status, all types of information data must be collected in the same vehicle, the same environment and the same time, the collection of equipment, and the type of object closely related to the type of detection. Therefore, the body of the vibration information data acquisition experiments and engine noise information data acquisition experiments were conducted in the same place, the two information data acquisition time points, and the same number of times. In the experiment, nine WX-C800 vibration sensors were set in the body of the bottom measurement and the left and right sides, each side of three, and evenly distributed. Engine noise information data collection time points and the same number of times. The measured data are listed in Table 6. Engine noise measurements were performed after 10 min of operation, and five samples were taken in the experiment, each with a time interval of 5 min. The measured data are presented in Table 6.

Table 6. Measured values of car body vibration information data (Hz)

| Sensor number | No. of measurements | | | | | Average value $F_j$ | Standard deviation $S_i$ |
|---|---|---|---|---|---|---|---|
| | NO.1 | NO.2 | NO.3 | NO.4 | NO.5 | | |
| C1# | 1.75 | 1.42 | 1.48 | 1.26 | 1.46 | 1.405 | 0.086 |
| C2# | 1.31 | 1.33 | 1.35 | 1.38 | 1.42 | 1.358 | 0.039 |
| C3# | 1.36 | 1.32 | 1.38 | 1.11 | 1.45 | 1.378 | 0.047 |
| C4# | 1.28 | 1.35 | 1.32 | 1.39 | 1.46 | 1.360 | 0.062 |
| C5# | 1.25 | 1.46 | 1.36 | 1.45 | 1.35 | 1.374 | 0.077 |
| C6# | 1.38 | 1.34 | 1.35 | 1.42 | 1.35 | 1.368 | 0.029 |
| C7# | 1.78 | 1.38 | 1.43 | 1.41 | 1.47 | 1.423 | 0.033 |
| C8# | 1.35 | 1.36 | 1.44 | 1.42 | 1.35 | 1.384 | 0.038 |
| C9# | 1.35 | 1.26 | 1.25 | 1.45 | 1.28 | 1.318 | 0.075 |
| Measured value $C$ | | | | | 1.374 | | |
| Standard deviation $S$ | | | | | 0.028 | | |

5.3.2 Information Data Analysis and Processing

In the processing of engine information data, we verified that the fractional-order differential algorithm improved the quality of the detection information data. Therefore, the fractional-order differential fusion algorithm was also applied to the processing of body vibration information data. Owing to the same testing location, working environment, and testing time point as the engine, the body vibration information data processing methods and processes with the engine information sound information data processing is the same, can be divided into the car sound vibration information data analysis and the car sound vibration information processing of the two parts of the article by the length of the article is now; this paper is only a simple description of the parts, no longer a detailed description of the content of the work of the various parts:

(1) Analysis of information data

In the experiment, we set the measurement error of the experimental data to 5%, i.e., the qualifying criterion of the data is 1.24<$C$<1.51. Therefore, data 1.75, 1.11 and 1.78 in Table 6 are obviously out of error and should be removed. A summary of the body vibration information data was obtained by analyzing the experimental data in Table 6.

Table 7. Summary of experimental data on car body vibration (dB)

| Sensor number | C1# | C2# | C3# | C4# | C5# | C6# | C7# | C8# | C9# |
|---|---|---|---|---|---|---|---|---|---|
| Average value $C_i$ | 1.405 | 1.358 | 1.378 | 1.360 | 1.374 | 1.368 | 1.423 | 1.384 | 1.374 |
| Standard deviation $S_i$ | 0.086 | 0.039 | 0.047 | 0.062 | 0.077 | 0.029 | 0.033 | 0.038 | 0.045 |
| measured value $C$ | | | | | 1.374 | | | | |

(2) Processing of information data
 ① Selection of parameters in the fractional-order differentiation algorithm: According to Table 7, choose $v=0.5$, $h=0.001$, then $n=[0.086-0.029]/0.001=57$
 ② Calculation of the relationship between the detection value and the impact factor: The application of the least squares method can be obtained between the bodywork detection data $C_i$ and the impact factor $S_i$ function formula:

$$C(x) = 37.0211*x^2 - 4.1877*x + 1.4844 \quad (19)$$

 ③ Establishment of an information data fusion model according to the Grunwald-Letnikov definition of fractional-order calculus. The body vibration information data fusion model based on the fractional-order differential algorithm can be obtained by combining the experimental data in Table 6 and (19):

$$C^v(x) \approx \lim_{h \to 0} h^{-v} \sum_{i=0}^{[(b-a)/h]} \binom{-v}{i} C(x-ih) = 0.01^{-0.5} \left[ C(x) - 0.5 C(x-0.01) - \frac{0.5(1-0.5)}{2!} C(x-2\times 0.5) - \ldots - \frac{0.5(1-0.5)(2-0.5)..(57-1-0.5)}{57!} C(x-57) \right] \quad (20)$$

 ④ By combining the fractional-order differentiation parameters in Table 7. and (13) and (14), the values of the bodywork experimental data processed by the 0.5 order differentiation algorithm can be obtained, as shown in Table 8. As you can see from it the accuracy is 4 times higher.

Table 8. Summary of experimental data of car body vibration after 0.5 order differential treatment (Hz)

| Sensor number | C1# | C2# | C3# | C4# | C5# | C6# | C7# | C8# | C9# |
|---|---|---|---|---|---|---|---|---|---|
| Post-fusion valuee $\overline{C_i^{0.5}}$ | 1.385 | 1.370 | 1.382 | 1.373 | 1.375 | 1.368 | 1.363 | 1.371 | 1.369 |
| Standard deviation after fusion $S^{0.5}$ | | | | | 0.007 | | | | |

5.3.3 Fault diagnosis and early warning

 Automobile body vibration during normal operation is a comprehensive embodiment of the operating state of the car; therefore, this paper will be the vibration measurement value of the body as the basis for judging the operating state of the body of the sound as the car noise measurement value of the iconic components of the faulty operating state of the car.
(1) Troubleshooting
 In this experiment, the experimental environment of the car body was the same as that of the engine, and the testing method and the number of tests of information data were in line with the international requirements for the testing conditions of the car body.
 According to the regulations of China for motor vehicles driving outside noise extreme value $C_b$ less than or equal to 1.5 Hz, because the entire car performance parameter of the early warning coefficient $k$ is 0.98, then $C_y = k*C_b = 1.47$. Referring to the accurate experimental data obtained in the experiment, the measured value $C = 1.374$, it can be determined that the current body is in a normal working condition.
(2) Fault warning
 Owing to the important parameter of acoustic prediction of engine failure, failure warning involves obtaining the point in time when the acoustic measurements are at the warning value set by the system. Therefore, the car measurement value $C$ is a function $C(t)$ at time t. The premise of accurate early warning of car failure is to accurately obtain the function $C(t)$, and then, according to the system, set warning acoustic measurements $C_y$ to derive the warning time point ty corresponding to it.
 In this study, we fit the function $C(t)$ according to the measured value Ji of the engine sound at different time points ti. It is reasonable that the more samples that are adopted, the more accurate the functional equation $C(t)$ is obtained, but this will lead to the shortcomings of high cost and low efficiency. Therefore, the samples used in the experiment were 6. The experimental samples were tested 0.5 months apart, and the experimental environment, information data collection, and processing methods were the same. The correspondence between the measured sample values $C_i$ and time is presented in Table 9.

Table 9. Measurements of automobile body vibration at different time points (Hz)

| Measured value $C_i$ | 1.358 | 1.362 | 1.376 | 1.383 | 1.395 | 1.412 |
|---|---|---|---|---|---|---|

| Time point $t_i$ (months) | 0 | 0.5 | 1.0 | 1.5 | 2.0 | 2.5 |
|---|---|---|---|---|---|---|

The least-squares method was again applied to fit the function $C(t)$ between the body vibration measurements $C_i$ and time $t_i$.

$$C(t)=0.0215*t+1.3541 \tag{21}$$

Combining the acoustic measurements $C_y = 1.47$ at the time of the engine fault warning and (21), we obtain

$$0.0215*t+1.3541=1.47 \tag{22}$$

Solving (22) shows that the warning time of body failure is $t_y=53.91$, i.e. about 53.9 months from the time of the first experiment. In other words, the warning time for the occurrence of failure was approximately 53.9 months according to the first measurement.

## 5.4 Automobile fault diagnosis and early warning

Automobiles consist of five major components: engines, gearboxes, bodies, chassis, and electrical systems. Under the provisions of the national standard, each component has a service life of the provisions of this life depending on the important performance parameters that lead to the failure of part of the change of law. Therefore, the normal operation process of the component is the process of continuous accumulation of performance parameters to the parameter to reach its rated value, which means that the component has been in a failure warning or failure state. Therefore, the automobile parts of Real-time operational status evaluation technology are the essence of the application of accurate information data and standard value comparison technology.

According to the previous discussion on the principle of automobile faults and early warning, as long as there is a component failure, it means that the entire car is in a state of failure. Therefore, the fault diagnosis and early warning process of the car involves predicting the time point of failure through the relevant information data of each component and the rated standard, and then comprehensively comparing the time point of failure of each component to obtain the most recent time point of the process.

Owing to limited space, this study only adopted the automobile engine and body as the standard of the five major components, and the abnormal sound and vibration as the main characteristics of the engine and body to identify its working status. Through precise calculations, it can be seen that under normal operating conditions (excluding emergencies), the warning of failure of the engine is 22.6 months after the first test time, and the warning of failure of the body is 53.9 months after the first test time. Therefore, it can be judged that the vehicle will have a failure warning 22.6 months after the first test time, and the type of failure warned will be an engine failure..

## 6. Conclusion

Cars have become the main means of transportation in people's lives, and in the process of driving, failure is directly related to the owner's safety. Therefore, it is of great significance to carry out fault warnings for automobiles to notify the owner in advance before the fault is about to occur, so as to give the owner sufficient time to maintain the vehicle. This paper proposes a method for real-time collection of relevant performance parameters of the main components of a car, predicting its failure warning time, and comparing the warning time of each component to obtain the first failure components, to realize the warning of automobile failure, and determine the type of upcoming failure.

Information data are the premise for an accurate warning of automobile failure. However, because of the complex structure of automobiles, the working environment is variable, which leads to unpredictable and irregular measurement errors in the information data measured by the automobile information data detection system, greatly affecting the reliability of the fault diagnosis and early warning system. In this study, based on the spectral characteristics of the fractional-order differential operator with enhanced signal strength and improved accuracy of information data, the fractional-order differential algorithm is applied to the fusion processing of automotive inspection information data, realizing the high-precision detection function of automotive inspection information data and providing technical support for improving the reliability of automotive fault diagnosis and prediction systems. Through the application of the method in automobile engines, body fault diagnosis, and early warning experiments, their fault diagnosis and early warning functions are realized. Through a comparison of their warning times, the diagnosis and warning function of the entire vehicle fault of the automobile was realized.

**Note:** As the research purpose of this paper is real-time operational status evaluation technology, only the engine and body were selected as representatives of the introduction and verification, and engineering practice can be based on the need to select the appropriate test object..

## Data availability statement

All data supporting the findings of this study are included in the article (and any supplementary files).

## Acknowlelgments

This project was supported by the National Natural Science Foundation, China (Nos. 51878005 and 51 778 004), and the Anhui Provincial Education Commission Foundation, China (No. KJ2020A0488).

## References.